\documentclass[12pt]{amsart}

\usepackage{amssymb}
\usepackage{latexsym}
\usepackage{verbatim}
\usepackage{fullpage}

\input {cyracc.def}
\font\tencyr=wncyr10 scaled \magstephalf
\font\tencyi=wncyi10 scaled\magstephalf
\font\tencysc=wncysc10 scaled\magstephalf
\def\rus{\tencyr\cyracc}
\def\rusi{\tencyi\cyracc}
\def\rusc{\tencysc\cyracc}

\newcommand{\reb}[1]{\mbox{\bf  (\ref{#1})}}
\newcommand{\re}[1]{\textrm  (\ref{#1})}

\renewenvironment{proof}
{\noindent {\sl Proof.}\quad }{\hfill
$\square$ \vskip1.1ex\noindent }

\newenvironment{proof*}
{\noindent {\sl Proof.}\quad }{\hfill
$\square$}

\renewcommand{\theequation}{\thesection .\arabic{equation}}
\renewcommand{\thesubsubsection}{\theequation .\arabic{subsubsection}}

\catcode`\@=\active
\catcode`\@=11
\def\@eqnnum{\hbox to
.01pt{}\rlap{\hskip-\displaywidth(\mathbf{\theequation})}}
\catcode`\@=12

\newenvironment{s}[1]
{ \vskip1.2ex \refstepcounter{equation}
\noindent {\bf \theequation\enspace #1.} \begin{sl}}{\end{sl}
\vskip1.1ex\noindent }

\newenvironment{rem}[1]
{ \vskip1.2ex \refstepcounter{equation}
\noindent {\bf \theequation\enspace {#1}.} }{ \vskip1.1ex\noindent }

\newcommand {\be}{{\frak b}}

\newcommand {\g}{{\frak g}}
\newcommand {\h}{{\frak h}}

\newcommand {\el}{{\frak l}}
\newcommand {\mm}{{\frak m}}
\newcommand {\n}{{\frak n}}

\newcommand {\p}{{\frak p}}
\newcommand {\q}{{\frak q}}
\newcommand {\rr}{{\frak r}}
\newcommand {\es}{{\frak s}}
\newcommand {\te}{{\frak t}}
\newcommand {\ut}{{\frak u}}

\newcommand {\gln}{{\frak gl}_n}
\newcommand {\gltn}{{\frak gl}_{2n}}
\newcommand {\sln}{{\frak sl}_n}
\newcommand {\glv}{{\frak gl}(V)}
\newcommand {\spv}{{\frak sp}(V)}
\newcommand {\spn}{{\frak sp}_{2n}}
\newcommand {\sov}{{\frak so}(V)}
\newcommand {\son}{{\frak so}_{n}}
\newcommand {\sono}{{\frak so}_{2n+1}}
\newcommand {\sone}{{\frak so}_{2n}}

\newcommand {\gD}{{\goth D}}

\newcommand {\gS}{{\goth S}}

\newcommand {\esi}{\varepsilon}
\newcommand {\ap}{\alpha}

\newcommand {\co}{{\mathcal O}}

\newcommand {\VV}{{\Bbb V}}
\newcommand {\ad}{{\mathrm{ad\,}}}

\newcommand {\Ann}{{\mathrm{Ann\,}}}

\newcommand {\ind}{{\mathrm{ind\,}}}
\newcommand {\Lie}{{\mathrm{Lie\,}}}

\newcommand {\rk}{{\mathrm{rk\,}}}

\newcommand {\trdeg}{{\mathrm{trdeg\,}}}
\newcommand {\tr}{{\mathrm{tr\,}}}

\newcommand {\GR}[2]{{\textrm{{\bf #1}}}_{#2}}

\newcommand {\un}{\underline}

\newcommand {\qu}{\hfill $\square$}
\newcommand {\beq}{\begin{equation}}
\newcommand {\eeq}{\end{equation}}

\renewcommand{\le}{\leqslant}
\renewcommand{\ge}{\geqslant}

\newfam\Bbbfam\newfam\eufam\newfam\eusfam%
\font\Bbbfont=msbm10 scaled 1200%
\font\olala=msam10 scaled 1200%
\font\frak=eufm10 scaled 1400%
\font\Bbbsmallfont=msbm8%
\textfont\Bbbfam=\Bbbfont\scriptfont\Bbbfam=\Bbbsmallfont
\font\euzw=eufm10 scaled 1200%
\font\euac=eufm7 scaled 1200%
\font\euacc=eufm7 scaled 1000%
\textfont\eufam=\euzw\scriptfont\eufam=\euac%
\scriptscriptfont\eufam=\euacc%
\font\euszw=eusm10 scaled 1200%
\font\eusac=eusm7 scaled 1200%
\font\eusacc=eusm7 scaled 1000%
\textfont\eusfam=\euszw\scriptfont\eusfam=\eusac%
\scriptscriptfont\eusfam=\eusacc%
\def\frak{\fam\eufam}%
\def\goth{\fam\eusfam}%
\def\Bbb{\fam\Bbbfam}%

\def\square{\hbox {\olala\char"03}}
\def\bbk{\hbox {\Bbbfont\char'174}}

\begin{document}
\setlength{\parskip}{2pt plus 4pt minus 0pt}
\hfill {\scriptsize September 1, 2004} 
\vskip1ex
\vskip1ex

\title[Seaweed subalgebras]{An extension of Ra\"\i s' theorem 
and seaweed subalgebras of simple Lie algebras }
\author[D.\,Panyushev]{\sc Dmitri I. Panyushev}
\thanks{This research was supported in part by R.F.B.I. Grant no. 
  02--01--01041 and CRDF Grant no. RM1--2543-MO-03}
\maketitle
\begin{center}
{\footnotesize
{\it Independent University of Moscow,
Bol'shoi Vlasevskii per. 11 \\
119002 Moscow, \quad Russia \\ e-mail}: {\tt panyush@mccme.ru }\\
}
\end{center}

\noindent
The ground field $\bbk$ is algebraically closed and of  characteristic zero.
Let $\q$ be a Lie algebra over $\bbk$ and $\xi\in\q^*$. Let
$\q_\xi$ denote the stationary subalgebra of $\xi$ in the coadjoint
representation of $\q$.
In other words, $\q_\xi=\{x\in\q\mid \xi([x,y])=0 \ \forall y\in \q\}$.
The {\it index\/} of $\q$, denoted $\ind \q$, is defined by
\[ 
   \ind\q =\min_{\xi\in\q^*}\dim\q_\xi \ .
\]
If $\q$ is an algebraic Lie algebra and 
$Q$ is an algebraic group with Lie algebra $\q$,
then $\ind\q$ equals the transcendence degree of 
the field of $Q$-invariant rational functions on $\q^*$.
If $\q$ is reductive, then $\q$ and $\q^\ast$ are isomorphic as 
$\q$-modules and hence $\ind\q=\rk\q$.  
It is an important invariant-theoretic problem to 
study index and, more generally, the coadjoint
representation for non-reductive  Lie algebras. 

A very interesting class of non-reductive  Lie algebras consists of
the so-called {\it seaweed subalgebras\/} of reductive Lie algebras.
This class includes all parabolic subalgebras, see Section~\ref{vodorosli} for the details.
In \cite{mmj}, we studied the index for seaweed subalgebras of 
classical simple Lie algebras. We 
obtained inductive formulae for the index that always apply 
in case of $\sln$ and $\spn$, and sometimes work for $\son$.
Then some complementary results in the orthogonal case
were obtained by Dvorsky \cite{dvor}. 
One of the goals of this paper is to show that inductive constructions from
\cite{mmj} allow, in fact, to obtain much stronger results.

We begin with a general result which concerns $\Bbb N$-graded Lie algebras
with at most three summands. Let $\h=\h(0)\oplus\h(1)\oplus\h(2)$ 
be such a Lie algebra. It is shown that under certain constraints
there is a subalgebra $\q\subset \h(0)$ such that 
{\em (a)} $\bbk(\h^*)^H$ is naturally isomorphic to $\bbk(\q^*)^Q$, and 
{\em (b)} if the action $(Q:\q^*)$ has a generic stabliser, then so does
$(H:\h^*)$, and these generic stabilisers are equal. 
Here $H$ and $Q$ are connected groups with Lie algebras $\h$ and $\q$, respectively.
This can be regarded as an extension of Ra\"\i s' theorem on the index 
of semi-direct products \cite{rais}.
(See Section~\ref{semidir} for the details.) Actually, an
$\Bbb N$-grading of a Lie algebra
sometimes allows us to prove that the coadjoint representation has no regular
invariants. And we prove that this is always the case for parabolic subalgebras
of semisimple Lie algebras.
This curious fact seems to have not been observed before. 

Using the result on gradings with at most three summands, 
we show that, for series $\sln$ and 
$\spn$, the coadjoint representation of any seaweed subalgebra
possesses some properties, similar to those of the adjoint 
representation in the reductive case. That is, 
if $\es$ is an arbitrary seaweed subalgebra of $\sln$ or $\spn$ with 
the corresponding connected group $S$,
then \\
\hbox to \textwidth{\refstepcounter{equation}{\bf (\theequation)}\label{0.1}%
\parbox{414pt}{
\begin{itemize} \item[\sf (i)] \ the field $\bbk(\es^*)^S$ is rational; 
\item[\sf (ii)] \  the representation $(S:\es^*)$  has a generic stabiliser whose
identity component is a torus.
\end{itemize} }
}
The proofs are based on the fact that, for the classical series,
any seaweed subalgebra admits a suitable $\Bbb N$-grading with at most three
summands.  For $\sln$ and $\spn$,
this grading always satisfies the necessary constraints, and we can argue 
by induction on $n$.
Unfortunately, this is not always the case for $\son$, so that we have only
partial results in the orthogonal case.
In fact, there is an example of a parabolic subalgebra of ${\frak so}_8$
such that its coadjoint representation has no generic stabiliser
\cite[Sect.\,3]{tayu1}. Our results for $\sln$ (resp. $\spn$ and $\son$)
are given in Section~\ref{gl} (resp. Section~\ref{sp}).
Actually, our theorem on 3-term gradings applies not only to classical Lie
algebras. One can present a number of other cases, where it works and yields 
the answer similar to Eq.~\reb{0.1}. A couple of examples of such sort is given
for the exceptional algebra of type $\GR{F}{4}$. Motivated by all these examples, 
we conjecture that the field of invariants of any seaweed subalgebra is always rational and
if a generic stabiliser exists, then its identity component is necessarily a torus.  

In Section~\ref{new}, we show that our general 
affirmative results for $\sln$ and $\spn$
are, in a sense, the best possible.
For any simple $\g$ such that the highest root, $\theta$, is fundamental,  
we give a uniform description of a parabolic subalgebra
such that its coadjoint representation has no generic stabiliser.
Let $\ap$ be the unique simple root that is not orthogonal to $\theta$.
Then we take $\p$ to be the minimal parabolic subalgebra corresponding to $\ap$.
Our proof is based on some curious relations between $\ap$ and the {\it canonical
string\/} of strongly orthogonal positive roots (alias: {\it Kostant's cascade 
construction}). 
For ${\frak so}_8$, we recover the above-mentioned example from \cite{tayu1}. 
Finally, we recall that the highest root is fundamental if and only
if $\g$ is neither $\sln$ nor $\spn$.
 
{\small
{\bf Acknowledgements.} Some results presented here were obtained during my visit
to the Institut Fourier (Grenoble) in February 2002. I would like thank Michel 
Brion for some interesting remarks after my talk there.
The paper was finished during my stay at the
Max-Planck-Institut f\"ur Mathematik (Bonn) in August 2004. 
}


\section{Invariant-theoretic preliminaries}
\label{prelim}

\noindent
Algebraic groups are denoted by capital Latin letters and their Lie
algebras are denoted by the corresponding lower-case Gothic letters.
\\[.6ex]
If an affine algebraic group $H$ acts regularly on an algebraic variety $X$, then
$H_x$ stands for the stabiliser of $x\in X$. Similarly, the stationary subalgebra of
$x$ in $\h=\Lie H$ is denoted by $\h_x$.
We say that the action $(H:X)$ {\it has a generic stabiliser\/}, 
if there exists
a dense open subset $\Omega\subset X$ such that all stabilisers $H_\xi$, $\xi\in \Omega$,
are conjugate in $H$. Then each of the subgroups $H_\xi$, $\xi\in\Omega$, 
is called a generic stabiliser.
Similarly, one defines the notion of a {\it generic stationary subalgebra}, 
which is a subalgebra of $\h$. Clearly, the existence of 
a generic stabiliser implies that of a generic stationary subalgebra.
That the converse is also true is proved by Richardson \cite[\S\,4]{Ri}.
The points in $\Omega$ are said to be {\it generic}.
The reader is referred to \cite[\S 7]{VP} for basic facts on generic stabilisers.

If two actions $(H_1:X_1)$ and $(H_2:X_2)$ are given, where
$H_1\subset H_2$ and $X_1\subset X_2$, then we say that their
generic stabilisers are equal, if {\sf (a)} both generic stabilisers exist and 
{\sf (b)} there exist
generic points $x_i\in X_i$ such that $(H_1)_{x_1}=(H_2)_{x_2}$.

Suppose now that $X$ is irreducible. Then $\bbk(X)^H$ stands for the field
of rational $H$-invariants on $X$. A celebrated theorem of {\sc M.\,Rosenlicht}
says that there is a dense open $H$-stable subset $\Omega\subset X$
such that $\bbk(X)^H$ separates the $H$-orbits in $\Omega$, see e.g.
\cite[2.3]{VP}.
In particular, $\trdeg \bbk(X)^H=\dim X-\max\dim_{x\in X} H{\cdot}x$.


\section{On the coadjoint representation
of some $\Bbb N$-graded Lie algebras.}  
\label{semidir}
\setcounter{equation}{0}

\noindent
Let $\h$ be an algebraic Lie algebra. Assume that it has an $\Bbb N$-grading of the form
$\h=\h(0)\oplus\h(1)\oplus\h(2)$. We also say that $\h$ has a {\it 3-term structure\/}.
Clearly, $\h(1)\oplus\h(2)$ is a nilpotent Lie algebra.
Therefore the algebraicity of $\h$ is equivalent to that of $\h(0)$.
Let $H(0)$ be a connected algebraic group with Lie algebra $\h(0)$.
Then $H=H(0)\rtimes \exp(\h(1)\oplus\h(2))$ is a connected group
with Lie algebra $\h$. The subspaces $\h(i)$ are $H(0)$-stable and
the decomposition $\h^*=\h(0)^*\oplus\h(1)^*\oplus\h(2)^*$ is therefore
$H(0)$-invariant, too. More precisely, the coadjoint representation of
$\h$, denoted $\ad^\ast$, satisfies the relation
$\ad^\ast(\h(i)){\cdot}\h(j)^*\subset \h(j-i)^*$.
Consider the following two conditions on $H$ and $\h$:
\begin{itemize}
\item[$(\Diamond_1)$] \ $H(0)$ has an open orbit in $\h(2)^*$, say $\co$;
\item[$(\Diamond_2)$] \ If $\xi\in\co$, then 
$H_\xi\subset H(0)\rtimes \exp(\h(2))$. In particular, 
$\h_\xi\subset\h(0)\oplus \h(2)$.
\end{itemize}
In the rest of the section, $\xi$ is an arbitrary but fixed point in $\co$.

\begin{s}{Theorem}   \label{main2}
Suppose $H$ satisfies conditions $(\Diamond_1)$ and
$(\Diamond_2)$. Then 
\begin{itemize}
\item[\sf (i)] \ There is a natural isomorphism
$\bbk(\h^*)^H\simeq \bbk(\h(0)^*_\xi)^{H(0)_\xi}$;
\item[\sf (ii)] \ Let $(f)$ be the ideal of the union of all divisors
in $\h(2)^*\setminus\co$. Then,  
regarding $f$ as a function on the whole of $\h^*$, we have
$ \bbk[\h^*]^H \subset \bbk[\h(0)^*_\xi]^{H(0)_\xi}\subset \bbk[\h^*_{(f)}]^H$.
In particular, if $\h(2)^*\setminus\co$ does not contain divisors, then
$\bbk[\h^*]^H\simeq \bbk[\h(0)^*_\xi]^{H(0)_\xi}$.
\item[\sf (iii)] \ If the action $(H(0)_\xi: \h(0)^*_\xi)$ has a generic stabiliser, 
then so does $(H:\h^*)$, and these generic
stablisers are equal.
\end{itemize}
\end{s}\begin{proof}
Our plan is to construct a section ${\goth S}\subset \h^*$
and a subgroup $\bar H\subset H$, acting on $\goth S$, such that 
$\bbk(\h^*)^H\simeq \bbk({\goth S})^{\bar H}\simeq \bbk(\h(0)^*_{\xi})^{H(0)_\xi}$ and
$\bbk[{\goth S}]^{\bar H}\simeq \bbk[\h(0)^*_{\xi}]^{H(0)_\xi}$. Recall that 
$\co=H(0){\cdot}\xi$.
\\[.5ex]
$1^o$. \ Set $\gS=\{\xi\}\times \h(0)^*$ and 
$\bar H=H(0)_\xi \rtimes \exp\h(2)\subset H$. 
Then $\gS$ is an affine subspace
in $\h^*$ and clearly it is $\bar H$-stable. From condition
$(\Diamond_2)$ it follows that 
\begin{gather}  \label{bir}
  H{\cdot}\eta\cap\gS=\bar H{\cdot}\eta  \\ 
\label{ravny}  H_\eta=\bar H_\eta
\end{gather}
for any $\eta\in\gS$. Furthermore,
\begin{equation}  \label{raznos}
 H{\cdot}\gS=\co \times (\h(1)^*\oplus \h(0)^*) \ .
\end{equation}
Since $H{\cdot}\gS$ is dense in $\h^*$,
it follows from Eq.~\re{bir} and Rosenlicht's theorem that
\begin{equation}   \label{fields}
\bbk(\h^*)^H\simeq \bbk(\gS)^{\bar H} \ .
\end{equation} 
$2^o$. \ Again, since $H{\cdot}\gS$ is dense in $\h^*$, the restriction homomorphism
$ \bbk[\h^*]^H \to \bbk[\gS]^{\bar H}$, 
$f\mapsto f\vert_\gS$, is injective.
Therefore $\bbk[\h^*]^H$ is identified with a subalgebra of $\bbk[\gS]^{\bar H}$. 
Notice that $H{\cdot}\gS$ is open in $\h^*$, and the complement $\h^*\setminus H{\cdot}\gS$ 
contains a divisor if and only if $\h(2)^*\setminus\co$ does.
If $g\in  \bbk[\gS]^{\bar H}$, then, in view of Eq.~\re{fields}, 
it extends to a rational function, say $\hat g$, on the whole of $\h^*$.
Since $H{\cdot}\gS$ is open and $\h^*$ is normal,  $\hat g$ may only have poles
on divisors in $\h^*\setminus H{\cdot}\gS$. From this we conclude that
\begin{itemize}
\item[] if $\h(2)^*\setminus \co$ 
contains divisors, and $(f)$ is the ideal of the union of all divisors,
then
\begin{equation}  \label{5}
    \bbk[\h^*]^H \subset \bbk[\gS]^{\bar H} \subset \bbk[\h^*_{(f)}]^H \ ;
\end{equation}
\item[] in particular, if $\h(2)^*\setminus \co$ does not
contain divisors, then 
\begin{equation}  \label{4}
    \bbk[\h^*]^H\simeq \bbk[\gS]^{\bar H} ;
\end{equation}
\end{itemize}
Here we extend $f$ to the whole of $\h^*$ using the natural projection 
$\h^*\to \h(2)^*$.
\\
$3^o$. \ Thus, we may forget about $H$ and $\h^*$ and work only
with the $\bar H$-action on $\gS$.
As $N:=\exp(\h(2)$ is a (commutative unipotent) normal subgroup of
$\bar H$, we first understand the structure of $N$-orbits in $\gS$.
It is easily seen that $N$ acts as a group of translation and all its
orbits in $\gS$ have one and the same dimension; i.e., for any
$(\xi,v)\in\gS$, we have 
\[
N{\cdot}(\xi,v)=(\xi, v+\ad^*(\h(2)){\cdot}\xi) \ .
\]
Here $\ad^*(\h(2)){\cdot}\xi)\subset \h(0)^*$ and the annihilator of
$\ad^*(\h(2)){\cdot}\xi$ in $\h(0)$ is $\h(0)_\xi$. 
Hence, all $N$-orbits
are parallel affine subspaces of dimension $\dim\h(0)-\dim\h(0)_\xi=
\dim\h(2)$. This implies that the mapping 
$\gS=\{\xi\}\times\h(0)^*\to \h(0)^*/\Ann(\ad^*(\h(2)){\cdot}\xi)
\simeq \h(0)^*_\xi$ is the geometric  quotient for 
the $N$-action on $\gS$. Thus,
\begin{equation}   \label{6}
 \bbk[\gS]^{\bar H}\simeq(\bbk[\gS]^N)^{H(0)_\xi}=\bbk[\h(0)^*_\xi]^{H(0)_\xi}
\quad \& \quad 
  \bbk(\gS)^{\bar H}\simeq (\bbk(\gS)^N)^{H(0)_\xi}=\bbk(\h(0)^*_\xi)^{H(0)_\xi}
\ .
\end{equation}
Now, combining Equations~\re{fields},\re{5},\re{4}, and \re{6} yields
parts (i) and (ii) in the Theorem.
\\
$4^o$. \ Suppose $Q\subset H(0)_\xi$ is a generic stabiliser for 
the $H(0)_\xi$-action on $\h(0)^*_\xi$. Then $Q$ is also a generic stabiliser
for the $\bar H$-action on $\gS$, since $N$ acts freely on
$\gS$.  Finally, it follows from Equations~\re{ravny} and \re{raznos}
that  generic stabilisers for $(\bar H:\gS)$ and $(H:\h^*)$ are equal.
\end{proof}%
{\bf Remarks.} 
1. \ In view of condition $(\Diamond_1)$, the polynomial
$f\in \bbk[\h(0)^*]$ cannot be $H(0)$-invariant. 
It is a semi-invariant of $H(0)$ with a non-trivial weight.
Its natural extension  to the whole of $\h^*$ is a semi-invariant of
$H$.

2. The group $H(0)_\xi$ can be disconnected.

3. In \cite[Prop.\,1.5]{mmj}, we obtained a formula for the index of algebras with 3-term 
structure satisfying condition~$(\Diamond_2)$ for generic points of $\h(2)^*$.  
The result was that $\ind\h=\ind \h(0)_\xi$.
Condition~$(\Diamond_1)$  was not considered there.
Here, having  stronger hypotheses, we proved a stronger result that the corresponding
fields of invariants are naturally isomorphic. A relationship between 
generic stabilisers is also new. Anyway, the point is that, 
for applications we have in mind, condition~$(\Diamond_1)$ is always satisfied.
\\[.7ex]
The simplest situation, where Theorem~\ref{main2} applies, is that
of semi-direct product. Let $Q$ be a connected algebraic group with Lie algebra
$\q$.
If $\rho: Q\to GL(V)$ is a finite-dimensional representation of $Q$,
then we denote the corresponding representation of $\q$ by the same letter.
The linear space $\q\times V$ has a natural structure of Lie algebra,
with bracket $[\ ,\ ]\widetilde{\ }$ defined by the equality
\[
[(s_1,v_1),(s_2,v_2)]\widetilde{\ }=([s_1,s_2], \rho(s_1)v_2-\rho(s_2)v_1)\ .
\]
The resulting Lie algebra is denoted by 
$\es=\q\rtimes_\rho V$. It is 
a semi-direct product of $\q$ and an Abelian ideal $V$.
It is a particular case of $\Bbb N$-graded Lie algebras considered above.
Namely, we have here $\h(0)=\q$, $\h(1)=0$, and $\h(2)=V$. Since condition~$(\Diamond_2)$
is trivially satisfied here, 
the following is a straightforward consequence of Theorem~\ref{main2}.

\begin{s}{Corollary} \label{obob}
Let $\es=\q\rtimes_\rho V$ be a semi-direct product as above. Suppose  
$Q$ has an open 
orbit in $V^*$, and $\xi\in V^*$ is a point in the open $Q$-orbit. Then
\begin{itemize}
\item[\sf (i)] \ 
$\bbk(\es^*)^S\simeq \bbk(\q^*_\xi)^{Q_\xi}$;
\item[\sf (ii)] \ 
If $V^*\setminus Q{\cdot}\xi$ does not contain divisors,
then
$\bbk[\es^*]^S\simeq \bbk[\q^*_\xi]^{Q_\xi}$.
\item[\sf (iii)] If the the coadjoint representation of $Q_\xi$ has a generic stabilizer, 
then the coadjoint representation of $S$ has, and these generic stabilisers are
equal. 
\end{itemize}
\end{s}%
There is a famous formula of M.\,Ra\"\i s for the index of semi-direct products
\cite{rais}.
When $Q$ has an open orbit in $V^*$, it amounts to the equality
$\ind\es=\ind \q_\xi$, i.e., 
$\trdeg \bbk(\es^*)^S=\trdeg\bbk(\q^*_\xi)^{Q_\xi}$.
Hence Corollary~\ref{obob}(i) can be regarded as a
refinement of Ra\"\i s' theorem in this situation.
\\[.7ex]
Let $\h$ be an $\Bbb N$-graded Lie algebra,
i.e., $\h=\displaystyle\bigoplus_{i=0}^d \h(i)$. We assume that $
\h_+:=\oplus_{i\ge 1}\h(i)\ne 0$. Then the group $H$ is a semi-direct
product of $H(0)$ and $H_+=\exp(\h_+)$.

\begin{s}{Lemma}   \label{inv=0}
Suppose $\h_+$ is a faithful $\h(0)$-module and there is a
(semisimple) element $x\in \h(0)$ such that $[x,y]=jy$ 
for any $j$ and $y\in\h(j)$. Then $\bbk[\h^*]^H=\bbk$.
\end{s}\begin{proof}
Let $T_1\subset H(0)$ be the one-dimensional torus with Lie algebra $\bbk x$.
Then $\bbk[\h^*]^H\subset \bbk[\h^*]^{T_1}= \bbk[\h(0)^*]$. 
Here $ \bbk[\h(0)^*]$ is regarded as subalgebra of $\bbk[\h^*]$ using the surjection of 
$H$-modules $\h^*\to \h(0)^*$.
On the other hand, no functions in $ \bbk[\h(0)^*]\setminus \bbk$ can
be $H_+$-invariant. Indeed, let 
$x_1,\ldots,x_m$ be a basis for $\h(0)$ and $F=F(x_1,\ldots,x_m)\in \bbk[\h(0)^*]$.
For any $y\in \h_+$ and $t\in\bbk$ we have 
\[
   \exp(ty)\cdot F=F +  t(y\ast F) + (\text{terms of higher degree in}\ t)\ .
\]

\noindent
Here $y\ast F=\displaystyle \sum_{i=1}^m \frac{\partial F}{\partial x_i}\,[y,x_i]$.  
Notice that $[y,x_i]\in \h_+$ and $\frac{\partial F}{\partial x_i}\in
\bbk[\h(0)^*]$. Hence,
the faithfulness guarantee us that for any $F$ there is $y$ such that $y\ast F\ne 0$.
Thus, \\ $\bbk[\h^*]^H\subset \bbk[\h^*]^{T_1}\cap \bbk[\h^*]^{H_+}=\bbk$.
\end{proof}%
It seems that the following interesting assertion was not noticed before.

\begin{s}{Corollary}  \label{invpar=0}
Let $\p$ be a proper parabolic subalgebra of a semisimple Lie algebra $\g$.
Then $\bbk[\p^*]^P=\bbk$.
\end{s}\begin{proof*} It is well known that any proper
parabolic subalgebra has a non-trivial $\Bbb N$-grading which is 
determined by a semisimple element in the centre of a Levi subalgebra of $\p$.
For this grading, the algebra $\p_+$ is just the nilpotent radical, $\p^{nil}$. 
\end{proof*}%
\begin{rem}{Example}    \label{primer-sl}
Let $\p$ be the maximal parabolic subalgebra of ${\frak gl}_{2n}$ whose
Levi subalgebra is isomorphic to $\gln\oplus\gln$.
In the matrix form, we have
$\p=\left\{\left(\begin{array}{cc} X & Y \\ 0 & Z \end{array}\right) \mid X,Y,Z\in\gln\right\}$.
The dual space $\p^*$ can be identified with $\gltn/\p^{nil}$, so that we
regard $\p^*$ as the set of matrices 
$\p^*=\left\{\nu=\left(\begin{array}{cc} \mm & \ast \\ \rr & \n \end{array}\right) 
\mid \mm,\rr,\n \in\gln\right\}$, where the contents of the right upper corner 
is irrelevant. We write a generic element of $P$ as
$p=\left(\begin{array}{cc} A & AB \\ 0 & D \end{array}\right)=
\left(\begin{array}{cc} A & 0 \\ 0 & D \end{array}\right)
\left(\begin{array}{cc} I_n & B \\ 0 & I_n \end{array}\right)$, 
where $A,D \in GL_n$, $B\in\gln$, and $I_n$ is the identity matrix.
The action of $p$ on $\p^*$ is given in the matrix form by the formulae
\[
   p: \ \left\{ \begin{array}{ccc}
     \mm & \mapsto & A\mm A^{-1}+A(B\rr)A^{-1} \\
     \n  & \mapsto & D\n D^{-1} - D(\rr B)D^{-1} \\
     \rr & \mapsto & D\rr A^{-1} \ .
\end{array}\right.
\]
It follows that $p$ takes $\rr\mm+\n\rr$ to
$D(\rr\mm+\n\rr)A^{-1}$. Therefore the matrix entries of $\rr$ and
$\rr\mm+\n\rr$ are regular invariants of the unipotent radical $P^u$,
i.e., the elements of the algebra $\Bbbk[\p^*]^{P^u}$.
Consider the open subset $\Omega\subset\p^*$, where $\rr$ is invertible.
On this open subset, $p$ takes the matrix
$\mm+\rr^{-1}\n\rr=\rr^{-1}(\rr\mm+\n\rr)$  to $A\bigl(\rr^{-1}(\rr\mm+\n\rr)\bigr)A^{-1}$.
Hence the rational functions 
$\hat g(\nu)=\tr(\mm+\rr^{-1}\n\rr)^i$, $i=1,2,\ldots,n$, lie in the field 
of invariants $\Bbbk(\p^*)^P$.
Using the scheme of the proof of Theorem~\ref{main2}, one
can prove that these functions generate the field
$\Bbbk(\p^*)^P$. Here $\p$ has a semi-direct product structure, i.e.,
$\p(0)=\gln\oplus\gln$, $\p(1)=0$, and $\p(2)=\p^{nil}$.
Then $\xi=\left(\begin{array}{cc} 0 & \ast \\ I_n & 0 \end{array}\right)$,
${\goth S}=
\left\{\left(\begin{array}{cc} \mm & \ast \\ I_n & \n \end{array}\right) 
\mid \mm,\n \in\gln\right\}$,
and $\bar{P}=\left\{\left(\begin{array}{cc} A & AB \\ 0 & A \end{array}\right)\mid
A\in GL_n,\ B\in \gln \right\}$. The algebra $\Bbbk[{\goth S}]^{\bar P}$ is freely generated
by the functions $g_i(\mm,\n)=\tr(\mm+\n)^i$, $i=1,\ldots,n$. 
It is easily checked that each $g_i$ extends to the rational
function $\hat g_i$  defined on $\Omega$.

Because in this example $\g$ is not semisimple, we have $\bbk[\p^*]^P\ne\bbk$.
Indeed, the algebra $\bbk[\p^*]^P$ is generated by
$\hat g_1(\nu)=\tr(\mm+\n)$.
This means that $\bbk[\p^*]^P=\Bbbk[{\goth S}]^{\bar P}$ if and only if
$n=1$, whereas the corresponding fields of invariants are always isomorphic.
Of course, the reason  is that the complement of $\Omega$ is a divisor.
\end{rem}%
%


\section{Seaweed subalgebras of simple Lie algebras}  
\label{vodorosli}
\setcounter{equation}{0}

\noindent
Seaweed subalgebras of $\gln$ are introduced in \cite{deki}. 
A general definition is given in \cite{mmj}.
Two parabolic subalgebras $\p$ and $\p'$
of $\g$ are said to be {\it weakly opposite\/}, if $\p+\p'=\g$.
Then the intersection $\es=\p\cap\p'$ is called a {\it seaweed subalgebra\/} 
of $\g$. If $\p$ and $\p'$ are opposite in the usual sense, then
$\es$ is a Levi subalgebra in either of them. At the other extreme,
if $\p'=\g$, then $\es=\p$.  
That is, theory of seaweed subalgebras can be regarded as a common generalisation
of the theory of parabolic subalgebras and Levi subalgebras.
 
Fix a Borel subalgebra $\be\subset\g$ and a Cartan subalgebra $\te$ in it, and
let $\be^-$ denote the opposite Borel subalgebra.
It was remarked in \cite{mmj} that any seaweed subalgebra is $G$-conjugate
to a subalgebra containing $\te$ and such that $\p\supset\be$ and
$\p'\supset\be^-$. Such seaweed subalgebras are said to be {\it standard\/}.
A standard seaweed subalgebra is determined by two subsets of the set
of simple roots of $\g$. It was conjectured in \cite{mmj} that $\ind\es\le \rk\g$ and
the equality holds if and only if $\es$ is a reductive (i.e., Levi) subalgebra.
The inequality is recently proved in \cite{tayu2}. 

In the next two sections, we consider seaweed subalgebras in classical Lie algebras. 
This can be regarded as a sequel to our article \cite{mmj}. The content of
these sections can be summarised in the following recipe: 

{\it Repeat the constructions
of \cite{mmj} and use Theorem~\ref{main2} (and Corollary~\ref{obob}) in place of
Proposition~1.5 in \cite{mmj} in order to obtain stronger conclusions.}


\section{Seaweed subalgebras of ${\frak gl}_n$}  
\label{gl}
\setcounter{equation}{0}

\noindent
It is harmless but technically easier to deal with $\gln$ in place of $\sln$.
Recall necessary results and notation from \cite{mmj}.
An ordered sequence of positive integers
$\un{a}=(a_1,\dots,a_m)$ is called a {\it composition\/} of the number
$\sum a_i$. The numbers $a_i$ are said to be the {\it parts\/} or
{\it coordinates\/} of the composition. 
Let $V$ be an $n$-dimensional $\bbk$-vector space.  
It is well known that there is a bijection
between the conjugacy classes of 
parabolic subalgebras of ${\frak gl}_n=\glv$ and 
the compositions of $n$. Under this bijection, 
the parabolic subalgebra that corresponds
to $(a_1,\dots,a_m)$ is one preserving a flag $\{0\}\subset
V_1\subset\dots \subset V_{m-1}\subset V_m=V$, where
$\dim V_i/V_{i-1}=a_i$. Then
a seaweed subalgebra of $\glv$ can be defined 
as the subalgebra preserving two
``opposite'' flags in $V$. 
\begin{rem}{Definition}  \label{seaweed:def}
Let $\un{a}=(a_1,\dots,a_m)$ and $\un{b}=(b_1,\dots,b_t)$ be two
compositions of $n=\dim V$. Fix a basis $(e_1,\dots,e_n)$ for $V$, and
consider two flags 
$\{0\}\subset V_1\subset\dots \subset V_{m-1}\subset V_m=V$ and
$V=W_0\supset W_1\supset\dots \supset W_{t-1}\supset W_t=\{0\}$,
where $V_i=\langle e_1,\dots, e_{a_1+\dots +a_i} \rangle$ and
$W_j=\langle e_{b_1+\dots+b_j+1},\dots, e_{n} \rangle$.
The subalgebra of $\glv$ preserving these two flags
is called a seaweed subalgebra of $\glv$ or a seaweed of degree
$n$. It will be denoted by $\es(\un{a}\mid \un{b})$. 
\end{rem}%
{\bf Remark.} A basis-free exposition requires an intrinsic definition
of ``opposite" flags. Two flags 
$\{0\}\subset V_1\subset\dots \subset V_{m-1}\subset V_m=V$ and
$V=W_0\supset W_1\supset\dots \supset W_{t-1}\supset W_t=\{0\}$
are called opposite, if
$\dim(V_i\cap W_j)=\max\{0, \dim V_i+\dim W_j -n\}$ for all $i,j$.
It is not hard to show that two flags are opposite if and only if 
there exists a basis for $V$ satisfying the properties of
Definition~\ref{seaweed:def}.
\\[1ex]
\begin{figure}[htb]
\begin{center}
\setlength{\unitlength}{0.02in}
\begin{picture}(110,103)(-5,0)
\multiput(0,0)(100,0){2}{\line(0,1){100}}
\multiput(0,0)(0,100){2}{\line(1,0){100}}
\qbezier[16](30,70),(40,70),(50,70)
\qbezier[24](30,70),(30,85),(30,100)
\multiput(30,70)(0,2){11}{\line(2,1){20}}
\multiput(34,70)(-4,22){2}{\line(2,1){16}}
\multiput(38,70)(-8,24){2}{\line(2,1){12}}
\multiput(42,70)(-12,26){2}{\line(2,1){8}}
%
\thicklines
\put(0,70){\line(1,0){30}}
\put(30,55){\line(1,0){15}}
\put(45,45){\line(1,0){10}}   
\put(55,20){\line(1,0){25}}  
\put(30,70){\line(0,-1){15}}   
\put(45,55){\line(0,-1){10}}
\put(55,45){\line(0,-1){25}}
\put(80,20){\line(0,-1){20}}

\put(50,100){\line(0,-1){50}}
\put(70,50){\line(0,-1){20}}
\put(50,50){\line(1,0){20}}
\put(70,30){\line(1,0){30}}
\put(0,70){\line(0,1){30}}
\put(100,0){\line(0,1){30}}
\put(0,100){\line(1,0){50}}
\put(80,0){\line(1,0){20}}

\thinlines
\put(-10,81){$a_1$}
\put(20,61){$a_2$}
\put(69,11){$a_m$}
\put(21,104){$b_1$}
\put(58,54){$b_2$}
\put(17,17){{\large $0$}}
\put(78,75){{\large $0$}}
\qbezier[75](100,0),(50,50),(0,100)
\end{picture}
\caption{A seaweed subalgebra of $\glv$} \label{seaweed:ris}
\end{center}
\end{figure}
\noindent
A standard seaweed algebra is depicted in Figure~\ref{seaweed:ris}.
It is convenient to think of seaweeds in $\gln$ as matrix algebras of such form. 
The following proposition immediately follows from the definition.
(A quick look on Figure~\ref{seaweed:ris} is also sufficient.)
\begin{s}{Proposition}  \label{(a,b)} 
\\ \indent
1. \ $\es(\un{a}\mid \un{b})$ is parabolic $\Leftrightarrow$ $\un{a}=(n)$
or $\un{b}=(n)$;  
\par
2. \ $\es(\un{a}\mid\un{b})$ is reductive $\Leftrightarrow$ $\un{a}=\un{b}$; 
\par
3. \ If $a_1+\dots+a_k=b_1+\dots+b_l$ for some $k< m$ and $l< t$, then
$\es(\un{a}\mid \un{b})$ is isomorphic to a direct sum of two proper
subalgebras,
either of which is a seaweed algebra (of smaller degree) in its
own sense. In particular,
if $a_1=b_1$, then $\es(\un{a}\mid \un{b})\simeq {\frak gl}_{a_1}\dotplus
\es(a_2,\dots,a_m\mid b_2,\dots,b_t)$.   \qu
\end{s}%
It was shown in \cite{mmj} that every seaweed subalgebra of $\glv$ has a 
semi-direct product structure satisfying the assumptions of Corollary~\ref{obob}.
We used that structure to derive inductive formulae for the index of seaweed 
subalgebras. That is to say, we kept track of only the transcendence degree of the 
field of invariants. Now, having at hand Theorem~\ref{main2}, we observe that the
very same procedure gives much more information. 

\begin{s}{Theorem}  \label{main-gl}
Let $\es$ be a seaweed subalgebra of $\gln$. Then
\begin{itemize}
\item[\sf (i)] \ $\bbk(\es^*)^S$ is a rational field;

\item[\sf (ii)] The action $(S:\es^*)$ has a generic stabiliser, which is a torus.
\end{itemize}
\end{s}\begin{proof}
Suppose that $\es=\es(\un{a}\mid\un{b})$, where $\un{a}=(a_1,\ldots,a_m)$ and
$\un{b}=(b_1,\ldots,b_t)$ are compositions of $n$. If  $a_1=b_1$, then we split up
$\es$ using Proposition~\ref{(a,b)}(3). Therefore, 
we may assume that $a_1<b_1$.
It was shown in \cite[Theorem\,4.2]{mmj} that $\es\simeq \q\rtimes_\rho V$, where
$\q=\es(\un{a}\mid a_1, b_1{-}a_1,b_2,\ldots,b_t)$ and $V$ is a commutative ideal 
of dimension $a_1(b_1-a_1)$. The ideal $V$ is represented by the stripped region in
figure~\ref{seaweed:ris}.
Here $Q$ has an open orbit  $\co\subset V^*$ and the 
stationary subalgebra $\q_\xi$, $\xi\in\co$, is isomorphic to 
\[  \begin{array}{lll}
- & \es(a_2,\dots,a_m\mid b_1-2a_1,a_1,b_2,\dots,b_t), & \ {\text{if\ }}
a_1\le b_1/2 ; \\
- & \es(2a_1-b_1,a_2,\dots,a_m\mid a_1,b_2,\dots,b_t), & \ {\text{if\ }}
a_1 > b_1/2  .
\end{array}
\]
It is also true that the stabiliser $Q_\xi$ is connected.
Therefore applying Corollary~\ref{obob} we reduce the problem to a seaweed of smaller size.
Eventually, we arrive at the case of a reductive seaweed subalgebra, where the assertion of 
the theorem is well-known to be true. 
\end{proof}%
However, this procedure does not always preserve the algebra of invariants.
So, it is not at all clear what can be said about $\bbk[\es^*]^S$.


\section{Seaweed subalgebras of ${\frak sp}_{2n}$ and $\son$}  
\label{sp}
\setcounter{equation}{0}

\noindent 
The results for symplectic and orthogonal Lie algebras
are quite similar up to a certain point. But after that one encounters 
with different phenomena. 
This is explained in \cite[Sect.\,5 \& 6]{mmj}. Here we recall
the main steps, but without reproducing all the notation and results.

We begin with an arbitrary seaweed subalgebra of $\g=\spv$ or $\sov$.

\begin{s}{Theorem}  \label{so+sp}
Let $\es$ be a seaweed subalgebra of $\spv$ (resp. $\sov$) that is not parabolic.
Then there exists a parabolic subalgebra $\p$ 
in ${\frak sp}(U)$ (resp. ${\frak so}(U)$) with $\dim U<\dim V$ such that 
\begin{itemize}
\item[\sf (i)] \  $\bbk(\es^*)^S$ is isomorphic to a purely transcendental
extension of\/ $\bbk(\p^*)^P$;
\item[\sf (ii)] \ If the action $(P:\p^*)$ has a generic stabiliser, then so does
$(S:\es^*)$, and these generic stabilisers are equal.
\end{itemize}
\end{s}\begin{proof}
We argue by induction on $\dim V$. We also say that $\dim V$ is the size of $\es$.
If $\es$ is a non-parabolic seaweed subalgebra of $\g$,
then there is an inductive procedure, similar to that described in Section~\ref{gl} for $\glv$, see \cite[Theorem\,5.2]{mmj}.
The inductive step replaces $\es$ with
a seaweed subalgebra of smaller size, say $\es'$. This $\es'$
may split into a direct
sum of a reductive subalgebra $\el$ and a seaweed
subalgebra $\q$ of an even smaller size. 
Since Corollary~\ref{obob} applies at each step, 
we conclude that {\em (a)}
$\bbk(\es^*)^S$ is a pure transcendental extension
of $\bbk(\q^*)^{Q}$ of degree equal to the rank of $\el$; {\em (b)}
if the action $(Q:\q^*)$ has a generic stabiliser, then so does
$(S:\es^*)$, and these generic stabilisers are equal.
If $\q$ is not parabolic, then we continue further with
$\q$. 
\end{proof}%
This result completely reduces the problem to considering parabolic subalgebras in 
$\spv$ and $\sov$. 
It is worth mentioning that $\dim V-\dim U$ is even, so that the above reduction
preserves the type ({\bf B} or {\bf D}) in the orthogonal case.

Next, any parabolic subalgebra of $\g$ has a suitable 3-term structure \cite{mmj}.
Hence, we may try to proceed further using results of 
Section~\ref{semidir} in full strength. 

$\bullet$ \ In the symplectic case, this 3-term structure always satisfies 
hypothesis $(\Diamond_1)$ and $(\Diamond_2)$,
see the proof of Theorem~5.5 in \cite{mmj}. 
So, a further reduction, using
Theorem~\ref{main2} is always possible. 
This leads to the following result.

\begin{s}{Theorem}  \label{main-sp}
Let $\p$ be a parabolic subalgebra of $\spn$. Then
\begin{itemize}
\item[\sf (i)] \ $\bbk(\p^*)^P$ is a rational field;

\item[\sf (ii)] The action $(P:\p^*)$ has a generic stabiliser whose identity
component is a torus.
\end{itemize}
\end{s}\begin{proof}
Let $\p$ be a parabolic subalgebra of $\g=\spn$ and 
$\p=\p(0)\oplus\p(1)\oplus\p(2)$ the 3-term structure introduced in 
Theorem~5.5 in \cite{mmj}. A brief description of it is as follows. 
Let $\ap_1,\ldots,\ap_n$ be the usual
set of simple roots for $\spn$, i.e., $\ap_i=\esi_i-\esi_{i+1}$, $i< n$, and
$\ap_n= 2\esi_n$. We may assume that $\p$ is standard. 
Let $r$ be the minimal index such that $\ap_r$ is not a root of the
standard Levi subalgebra of $\p$. 

If $r=n$, then $\p$ is a maximal parabolic subalgebra, and we obtain 
the symplectic analogue of Example~\ref{primer-sl}. Here everything can be computed
quite concretely, see Example~\ref{primer-sp} below.

If $r < n$, we define the $\Bbb Z$-grading of $\g$ by letting $\g(i)$ be the
the sum of all root spaces $\g_\gamma$ with $[\gamma:\ap_r]=i$.
Here  $[\gamma:\ap_r]$ is the coefficient of $\ap_r$ in the expansion
of $\gamma$ via the simple roots. (Of course, the Cartan subalgebra $\te$ is included in 
$\g(0)$.) Restricting this grading to $\p$, we obtain the required $\Bbb N$-grading.
It follows from the construction that $\p(0)\simeq {\frak gl}_r\oplus \p'$,
where $\p'$ is a parabolic subalgebra in ${\frak sp}_{2n-2r}$.
Applying Theorem~\ref{main2}, we reduce the problem to the subalgebra
$\p(0)_\xi$, where $\xi\in\p(2)^*$ is a point in the dense $P(0)$-orbit.
Here $\p'$ acts trivially on $\p(2)$, and, as ${\frak gl}_r$-module,
$\p(2)$ is isomorphic to the space of symmetric $r\times r$-matrices.
It follows that 
$\p(0)_\xi\simeq {\frak so}_{r}\oplus \p'$. 
Here $P(0)_\xi=O_r\times P'$, i.e., it is disconnected. However disconnectedness
of $O_r$ only results in the fact that 
a generic stabiliser for $(P(0)_\xi:\p(0)_\xi)$, and hence for $(P:\p^*)$,
appears to be disconnected, but this
does not affect further reduction steps applied to $\p'$.
\end{proof}%
Combining Theorems~\ref{so+sp} and \ref{main-sp}, we obtain a complete 
answer for the seaweed subalgebras in $\spv$.

\begin{rem}{Example}  \label{primer-sp}
Let $\p$ the maximal parabolic subalgebra of $\spn$ whose Levi subalgebra is 
isomorphic to $\gln$. Choose a basis for $V=\bbk^{2n}$ such that the 
skew-symmetric bilinear form has the matrix
$\left( 
  \begin{array}{cl} 0 & -I_n \\ I_n & \hphantom{-}0 \end{array}
\right)$. Then 
$\p=\left\{\left(\begin{array}{cc} X & Y \\ 0 & -X^t \end{array}\right) 
\mid X,Y\in\gln,\ Y=Y^t \right\}$.
The dual space $\p^*$ can be identified with $\spn/\p^{nil}$, so that we
regard $\p^*$ as the set of matrices 
$\p^*=\left\{\nu=\left(\begin{array}{cc} \mm & \ast \\ \rr & -\mm^t \end{array}\right) 
\mid \mm,\rr \in\gln, \, \rr=\rr^t
\right\}$, where the contents of the right upper corner 
is irrelevant. We write a generic element of $P$ as
$p=\left(\begin{array}{cc} A & AB \\ 0 & (A^t)^{-1} \end{array}\right)$,
where $A \in GL_n$ and $B=B^t$.
The action of $p$ on $\p^*$ is given in the matrix form by the formulae
\[
   p: \ \left\{ \begin{array}{ccc}
     \mm & \mapsto & A\mm A^{-1}+A(B\rr)A^{-1} \\
     \rr & \mapsto & (A^t)^{-1}\rr A^{-1} \ .
\end{array}\right.
\]
It follows that $p$ takes $\rr\mm-\mm^t\rr$ to
$(A^t)^{-1}(\rr\mm-\mm^t\rr)A^{-1}$. Therefore the matrix entries of $\rr$ and
$\rr\mm-\mm^t\rr$ are regular invariants of the unipotent radical $P^u$,
i.e., the elements of the algebra $\Bbbk[\p^*]^{P^u}$.
Consider the open subset $\Omega\subset\p^*$, where $\rr$ is invertible.
On this open subset, $p$ takes the matrix
$\mm-\rr^{-1}\mm^t\rr$
to $A(\mm-\rr^{-1}\mm^t\rr)A^{-1}$.
Hence the rational functions 
$\hat g(\nu)=\tr(\mm-\rr^{-1}\mm^t\rr)^{2i}$, $i=1,2,\ldots,[n/2]$, lie in the field 
of invariants $\Bbbk(\p^*)^P$. (Clearly, the trace of an odd power equals zero.) 

Using the scheme of the proof of Theorem~\ref{main2}, one
can prove that these functions generate the field
$\Bbbk(\p^*)^P$. 
Here $\p$ has a semi-direct product structure, i.e.,
$\p(0)=\gln$, $\p(1)=0$, and $\p(2)=\p^{nil}$.
Then 
$\xi=\left(\begin{array}{cc} 0 & \ast \\ I_n & 0 \end{array}\right)$, 
${\goth S}=
\left\{\left(\begin{array}{cc} \mm & \ast \\ I_n & -\mm^t \end{array}\right) 
\mid \mm \in\gln\right\}$
and $\bar{P}=\left\{\left(\begin{array}{cc} A & AB \\ 0 & A \end{array}\right)\mid
A\in O_n,\ B=B^t \right\}$. The algebra $\bbk[{\goth S}]^{\bar P}$ is freely generated
by the functions $g_i(\mm)=\tr(\mm-\mm^t)^{2i}$, $i=1,\ldots,[n/2]$.
(For, the mapping $\mm\mapsto (\mm-\mm^t)$ is the factorisation with respect to
the action of $\bar P^u$ on $\gS$, and then one has to take $O_n$-invariants of the quotient
obtained.) 
It is easily checked that each $g_i$ extends to the rational
function $\hat g_i$  defined on $\Omega$.

Since $P(0)_\xi=O_n$, the identity component of the generic stabiliser is a torus
(of dimension $[n/2]$).
\end{rem}%
$\bullet$ \ In the orthogonal case, the 3-term structure that we constructed in
\cite{mmj} satisfies $(\Diamond_1)$, but does not
always satisfy $(\Diamond_2)$. 
Let $\ap_1,\ldots,\ap_n$ be the usual
set of simple roots for either $\sono$ or $\sone$, i.e., 
$\ap_i=\esi_i-\esi_{i+1}$, $i< n$, and
$\ap_n= \left\{ \begin{array}{rl} \esi_{n-1}+\esi_n & \text{ if } \ \g=\sone \\
                                   \esi_n & \text{ if } \ \g=\sono
\end{array}\right.$. 
\\[.6ex]
Let $\p$ be a standard parabolic subalgebra of ${\frak so}_N$, $N=2n$ or $2n+1$,
and let
$r$ be the minimal index such that $\ap_r$ is not a root of the 
standard Levi subalgebra of $\p$. 
Using this root, one constructs an $\Bbb N$-grading of $\p$, as above.
Then  condition 
$(\Diamond_2)$ is satisfied for this $\Bbb N$-grading if and only if $r$ is even,
modulo the following adjustment for $\sone$. If $n$ is even, then both $r=n-1,\, n$
are acceptable; if $n$ is odd, then neither of them is acceptable.
This yields the following  assertion.

\begin{s}{Theorem}  \label{main-so}
Let $\p$ be a standard parabolic subalgebra in ${\frak so}_N$. 
Let $\ap_{r_1},\ap_{r_2},\ldots $ be all simple roots that are not in
the standard Levi subalgebra of $\p$. If all numbers $r_1,r_2,\ldots$
are even (modulo the above adjustment for $\sone$), then
\begin{itemize}
\item[\sf (i)] \ $\bbk(\p^*)^P$ is a rational field;

\item[\sf (ii)] The action $(P:\p^*)$ has a generic stabiliser, which is a torus.
\end{itemize}
\end{s}\begin{proof} 
We argue as in the proof of Theorem~\ref{main-sp}, starting with $r=r_1$.
The inductive step bring us from 
$\p(0)={\frak gl}_r\oplus \p'$ to $\p(0)_\xi={\frak sp}_r\oplus \p'$.
Here $\p'$ is a parabolic subalgebra of ${\frak so}_{N-2r}$ and 
$P(0)_\xi=Sp_{r}\times P'$ is connected. Therefore
a generic stabiliser at the very end will be connected, too. It is also clear from this
outline, why it is important that $r$ is even.
\end{proof}%
Thus, our inductive method does not apply to every parabolic 
subalgebra of $\sov$.
Recent results of Tauvel and Yu \cite{tayu1} show that this is not 
a drawback of our approach. 
For, they constructed an example of a parabolic subalgebra in ${\frak so}_8$
such that the coadjoint representation does not have generic stabilisers.
Namely, $\p$ is the minimal parabolic subalgebra corresponding to $\ap_2$
(the branching node
on the Dynkin diagram). Here $\dim\p=17$ and the maximal dimension of $P$-orbits in 
$\p^*$ is 16. That is, there is a dense open subsets of $\p^*$ consisting of
a 1-parameter family of $P$-orbits of dimension 16.
The stabiliser of each orbit is a 1-dimensional unipotent  subgroup.
But these subgroups are not conjugate in $P$.
Still, the field of invariants in this example is rational, in view of 
the L\"uroth theorem.
Moreover, no examples is known with a non-rational field of invariants for the
coadjoint representation of a parabolic (or seaweed) subalgebra.
\\[.7ex]
One can notice that whenever our inductive procedure applies,
it gives the rationality of the field of invariants and the existence 
of a generic stabiliser. Furthermore, 
the identity component of a generic stabiliser appears to be a torus.

\begin{rem}{Example}   \label{osobye}
Let $\g$ be an algebra of type $\GR{F}{4}$.
Take a maximal parabolic subalgebra whose Levi subalgebra is 
of semisimple type $\GR{B}{3}$ or $\GR{C}{3}$. Then 
the natural $\Bbb N$-grading of both this parabolics satisfies Theorem~\ref{main2}.
This is because both ($2\,0{\Rightarrow} 0\,0$) and ($0\,0{\Rightarrow} 0\,2$)
are weighted Dynkin diagrams of quadratic nilpotent elements.
("quadratic" means that $(\ad e)^3=0$, see \cite[Example\,1.6]{mmj} about this.)
Here $\p(0)_\xi$ is
isomorphic to 
$\left\{\begin{array}{cc} {\frak so}_6 &  \text{ in the $\GR{B}{3}$-case} \\
                          {\frak sp}_6 &  \text{ in the $\GR{C}{3}$-case}
\end{array}\right.$.
From this we immediately obtain that in both cases
the field of invariants is rational,
and  generic stabilisers are three-dimensional tori. 

\end{rem}%
Based on these observations, we propose the following 

\begin{s}{Conjecture}
Let $\g$ be a simple Lie algebra and $\es\subset\g$ a seaweed subalgebra.
Then 
\begin{itemize}
\item[\sf (i)] \ the field\/ $\bbk(\es^*)^S$ is rational;
\item[\sf (ii)] \ If a generic stabiliser for $(S:\es^*)$ exists, then
its identity component is a torus.
\end{itemize}
\end{s}%
%


\section{Constructing coadjoint representations without generic stabilisers}  
\label{new}
\setcounter{equation}{0}

\noindent
In this section, $\g$ is a simple Lie algebra with a fixed triangular 
decomposition $\g=\ut^+\oplus\te\oplus\ut^-$. The corresponding set of 
roots (resp. positive roots) is $\Delta$ (resp. 
$\Delta^+$), and the highest root is $\theta$.
If $\gamma\in \Delta$, then $\g_\gamma$ is the corresponding root space  and 
$e(\gamma)$ is a nonzero vector in $\g_\gamma$.
As usual, we assume that all roots live in a $\Bbb Q$-vector space $\VV$
of dimension $\rk \g$,
and that $\VV$ is equipped with a $W$-invariant scalar product
$(\ ,\ )$, where $W$ is the Weyl group.

Recall the construction of the {\it canonical string\/} of 
(strongly orthogonal) roots in $\Delta^+$.
Sometimes, it is called {\it Kostant's cascade construction\/}.
We start with $\mu_1=\theta$, and then consider 
$\Delta_1=\{\gamma\in \Delta \mid (\gamma,\mu_1)=0\}$.
Here $\Delta_1$ is a root system in its own right, which can be reducible.
If $\Delta_1=\sqcup_{j=2}^s \Delta_1^{(j)}$, then we choose the highest root in
each irreducible subsystem. These are the following elements of the canonical string:
$\mu_2,\ldots,\mu_{s}$; here $\mu_{2}\in \Delta_1^{(2)}$, etc.
Then we do the same thing with each 
$\Delta_1^{(j)}$, and so on\dots
This procedure eventually terminates, and 
we obtain the canonical string ${\goth D}=\{\mu_1,\ldots,\mu_l\}$.
Each member of $\goth D$ is the highest root in a certain irreducible subsystem
of $\Delta$, and the roots in $\gD$ are pairwise mutually strongly orthogonal.
In particular, the roots in $\gD$ are linearly independent
and $l\le \rk\g$.
 It is clear that the numbering of roots arising in each step is not essential.
What is essential is a poset structure of on $\gD$. Namely, 
$\mu_1$ is the unique maximal element, and the elements covered by $\mu_1$
are precisely $\mu_2,\dots,\mu_{s}$. The elements of $\gD$ covered by  
$\mu_2$ are precisely the highest roots of the irreducible subsystems
of $\{\gamma\in \Delta_1^{(2)}\mid (\gamma,\mu_2)=0\}$, and likewise in each step.
The Hasse diagrams of these posets for all $\g$ can be found in \cite[Table~III]{joseph}.

The canonical strings are of interest for us because of the following 
result. We may identify the dual space $\be^*\simeq \g/\ut^+$ with 
the vector space $\be^-=\ut^- \oplus\te$. In order to distinguish the true
$\be^-$ and $\be^*$, a nonzero root vector in $\be^*$ corresponding to
a root $\gamma$ is denoted by $\xi(\gamma)$.

\begin{s}{Proposition}   \label{borel}
The vector  
$\displaystyle
\xi_0=\sum_{i=1}^l \xi(-\mu_i)
$
is a generic point in $\be^*$ and the identity component of
$B_{\xi_0}$ is a torus of dimension\/ $\rk\g- l$. Actually,
$\be_{\xi_0}$ is equal to 
$\h=\{x\in \te\mid \mu_i(x)=0 \ \ \forall i\}$.
Furthermore, if we regard $\h$ as a subspace of $\be^*$, 
then each point in the affine subspace $\xi_0+\h$ is generic and
$B{\cdot}(\xi_0+\h)$ contains a dense open subset of $\be^*$.
\end{s}%
This is implicit in Joseph's article \cite{joseph}, and was known for many years
as a folklore. In fact, Joseph shows that $\xi$ is a generic point in the
$\be$-module $(\ut^+)^*\simeq \g/\be$. Then a minor adaptation of his arguments,
together with rudiments of invariant-theoretic technique, is sufficient 
to get the above result. In case of $\gln$, a proof of Proposition was given in 
\cite[\S\,3]{arh}.
A general proof based on the cascade construction
is given in \cite[Theorem\,3.7]{GS}. 
Essentially the same proof appears recently in \cite{tayu2}.


\begin{s}{Lemma}     \label{w_0}
Let $w_0$ be the longest element of $W$. 
Consider the subspace $\VV'=\{ x\in \VV\mid  w_0(x)=-x\}$.
Then the elements of $\gD$ form a basis for $\VV'$.
\end{s}\begin{proof}
It can be shown {\sl a priori\/} (or by a straightforward verification)
that $l=\dim\VV'$. In particular, $l=\rk\g$ if and only if $\VV'=\VV$.
Therefore, it remains only to verify that for $\g\in\{
\GR{A}{n},\,\GR{D}{2n+1},\, \GR{E}{6}\}$ the roots $\mu_i$ lie in $\VV'$.
\end{proof}%
Now, we are ready to provide a series of examples.

\begin{s}{Theorem}  \label{main6}
Suppose $\g$ is such that $\theta$ is a fundamental weight.
Let $\ap$ be the unique simple root that is not orthogonal to \,$\theta$.
Let $\p=\be\oplus \g_{-\ap}$ be the standard minimal parabolic subalgebra 
corresponding to $\ap$.
Then $\ind\p=\ind\be+1$ and the coadjoint representation $(P:\p^*)$
does not have a generic stabiliser.
\end{s}%
Practically, the hypothesis on fundamentality 
means that $\g$ is neither $\sln$ nor $\spn$. Therefore, taking into account our results
in Sections~\ref{gl} and \ref{sp}, we obtain

\begin{s}{Corollary}
Given a simple Lie algebra $\g$, the following conditions are equivalent:
\begin{itemize}
\item[\sf (i)] \ $\g=\sln$ or $\spn$;
\item[\sf (ii)] \  For any seaweed subalgebra $\es\subset\g$, the
coadjoint representation $(S:\es^*)$ has a generic stabiliser.
\end{itemize}
\end{s}%
{\sl Proof of Theorem~\ref{main6}.}
\quad
The proof occupies the rest of this section. It exploits an interesting
relation between $\gD$ and $\ap$.

Recall that $(-w_0)$ is an involutory linear transformation of $\VV$,
and $(-w_0)\Delta^+=\Delta^+$. Since $(-w_0)\theta=\theta$, it follows from our 
hypothesis that $(-w_0)\ap=\ap$, as well. 
Therefore,  by Lemma~\ref{w_0}, $\ap$ lies in the\/ $\Bbb Q$-linear span of $\gD$.
Hence,  $\ap=\sum_{i}k_i\mu_i$. We are interested in the coefficients of this
expansion.

\begin{s}{Lemma}  \label{koeff}
\begin{itemize}
\item[\sf 1.] \  $k_1=1/2$; \ 
 if $k_i\ne 0$ for $i\ge 2$, then $(\mu_i,\ap)<0$.
\item[\sf 2.] \ If $i\ge 2$ and $(\mu_i,\ap)<0$, then $k_i
=-\displaystyle\frac{\|\ap\|}{2\|\mu_i\|}$ is negative.
\item[\sf 3.] \ $\sum_{i\ge 2} k_i=-\frac{3}{2}$ or, equivalently,
$\sum_{i} k_i=-1$.
\end{itemize}
\end{s}\begin{proof*} 
1. Recall that $\mu_1=\theta$.
By the assumption, we have
$1=(\ap, \theta^\vee)=k_1(\theta,\theta^\vee)=2k_1$.
Since the $\mu_i$'s are pairwise mutually orthogonal, $k_i\ne 0$ 
if and only if $(\ap,\mu_i)\ne 0$. 
If $(\ap,\mu_i)\ne 0$, $i\ge 2$, then this number cannot be positive. For, otherwise
$\mu_i-\ap$ would be a positive root and then $0\le (\theta,\mu_i-\ap)=
-(\theta,\ap)<0$, which is absurd.

2.  Since $\theta$ is fundamental, we have $(\theta,\ap^\vee)=1$.
That is, $\ap$ is necessarily a long root. Therefore,
if $(\mu_i,\ap)<0$, then actually, $(\mu_i,\ap^\vee)=-1$.
Hence 
\[
2k_i=k_i(\mu_i,\mu_i^\vee)=(\ap,\mu_i^\vee)=-\displaystyle\frac{\|\ap\|}{\|\mu_i\|} .
\]
3. Now,
\[
  2=(\ap,\ap^\vee)=k_1(\mu_1,\ap^\vee)+\sum_{i\ge 2}k_i(\mu_i,\ap^\vee)=
    1/2- \sum_{i\ge 2}k_i \ ,
\]
and we are done.
\end{proof*}%
\begin{s}{Corollary}   \label{I}
Set $I=\{ i \mid (\mu_i,\ap)<0\} \subset \{2,\ldots,l\}$. Then $\#I\le 3$
and 
\[
   \ap= \frac{1}{2}(\theta- \sum_{i\in I} \frac{\|\ap\|}{\|\mu_i\|}\mu_i)=
\sum_{i\in I\cup\{1\}}k_i\mu_i \ .
\]
\end{s}%
{\bf I.} \  For convenience, we first consider the case in which 
$\ind\be=0$. This means that $l=\rk\g$ and $\h$, the space introduced in 
Proposition~\ref{borel}, is zero. 
Although this is not needed for our proof, we notice that this means
that $\g\in\{ \GR{B}{n}\ (n\ge 3)
,\,\GR{D}{2n}\,(n\ge 2),\,\GR{E}{7},\,\GR{E}{8},\,\GR{F}{4},\,\GR{G}{2}\}$.

Again, we identify the dual space $\p^*=\g/\p^{nil}$
with the space $\be^-\oplus \g_\ap \subset \g$, with the same notation
for root vectors in $\p$ and $\p^*$. 
\\[.6ex]
Set $\xi_a=\sum_{i=1}^l \xi(-\mu_i) + a \xi(\ap)$, $a\in\bbk$.
We are going to prove that the affine line
$L=\{\xi_a \mid a\in\bbk\}$ has the property that
$P{\cdot}L$ contains a dense open subset of $\p^*$; 
$\dim\p_{\xi}=1$ for every $\xi\in L$ and  
neither of the stabilisers $\p_\xi$ can be a generic stabiliser.

Notice that the image of $\xi_a$ in $\be^*$ is the generic point given in
Prop.~\ref{borel}. Therefore $\dim P{\cdot}\xi_a\ge \dim B$ for all $a$.
Since $\dim\p=\dim\be+1$,  we conclude, 
for parity reasons, that $\dim P{\cdot}\xi_a= \dim B$ 
and $\dim \p_{\xi_a}=1$ for all $a$.
We are going to give explicit expressions for all these stationary subalgebras.
To simplify the formulae, we assume that the root vectors in $\p^*$ are
already fixed, but the explicit choice (normalisation)
of vectors $e(\gamma)\in\p$ is still at our disposal. We use the notation of 
Corollary~\ref{I}.

\begin{s}{Proposition}   \label{stab}
Under a suitable choice of root vectors, 
the one-dimensional space $\p_{\xi_a}$ is generated by 
\[
e_a=e(-\ap)+ \sum_{i\in I}e(\theta{-}\ap{-}\mu_i)- a\sum_{i\in I\cup\{1\}}e(\mu_i)= 
\sum_{i\in I\cup\{1\}}e(\theta{-}\ap{-}\mu_i)- a\sum_{i\in I\cup\{1\}}e(\mu_i) \ .
\]
\vskip-1ex
\end{s}\begin{proof*}
1. We begin with the case $a=0$.
Computing the expression $e_0{\cdot}\xi_0$ we obtain
\begin{multline*}   
 \Bigl(e(-\ap)+ \sum_{i\in I}e(\theta-\ap-\mu_i)\Bigr){\cdot} 
\Bigl(\sum_{i=1}^l \xi(-\mu_i)\Bigr)= \\
=\sum_{i\in I}
\Bigl( e(-\ap){\cdot}\xi(-\mu_i)+ e(\theta-\ap-\mu_i){\cdot}\xi(-\theta) \Bigr)
+ \Bigl(\sum_{i\in I} e(\theta-\ap-\mu_i)\Bigr){\cdot}
\Bigl( \sum_{i=2}^l \xi(-\mu_i)\Bigr) \ .
\end{multline*}
In the passage to the second row we used the fact that $e(-\ap)\cdot\xi(-\mu_i)\ne 0$
if and only if $\ap+\mu_i$ is a root, i.e., $i\in I$.
It is clear that under suitable choice of $e(\gamma)$'s each summand of
the first sum in the second row can be made zero. As for the second sum, it is just equals
zero. To see this, we show that $\eta_{ij}:=\theta-\ap-\mu_i-\mu_j$  never belongs
to $\Delta^-\cup\{0\}\cup\{\ap\}$. Indeed,
$(\eta_{ij},\theta^\vee)=1$. Hence $\eta_{ij}\ne 0$, and if it is a root, then it must be
$\ap$. But $(\eta_{ij},\ap^\vee)=1-2+1-(\mu_j,\ap)\ne 2$, since
$\mu_j\ne \ap$. Hence $\eta_{ij}\ne\ap$.

2. Now, we consider $\xi_a$ with an arbitrary $a\in \bbk$.
The root vectors $e(\theta-\ap-\mu_i)$, $i\in I\cup\{1\}$,
are already chosen, but all other are still at our disposal.
Computing the expression $e_a{\cdot}\xi_a$ and using the fact that $e_a{\cdot}\xi_a=0$, 
we obtain
\begin{multline}   \label{xi_a} 
\bigl(e_0-a\sum_{i\in I\cup\{1\}}k_i e(\mu_i) \bigr){\cdot} 
\bigl(\xi_0 +a\xi(\ap)\bigr)= 
a e_0{\cdot}\xi(\ap)- a \sum_{i\in I\cup\{1\}} e(\mu_i){\cdot}\xi_0 =\\
=a e(-\ap){\cdot}\xi(\ap) -a\sum_{i\in I\cup\{1\}} e(\mu_i){\cdot}\xi(-\mu_i) \ .
\end{multline}
It is easily seen that all other summands are equal to zero.
For instance, $e(\mu_i){\cdot}\xi(\ap)=0$, since 
$\ap+\mu_i$ is either not a root at all, or not a root
of $\p^*$.
Also,  for $i\in I$ we have $e(\theta-\ap-\mu_i){\cdot}\xi(\ap)=0$, because 
$\theta-\mu_i$ is not a root.
\\[.6ex]
Now, the last expression in Eq.~\re{xi_a} is a sum of elements lying in
$\te\in \p^*$. 
Under the identification of $\te$ and $\te^*$, we have
$e(-\ap){\cdot}\xi(\ap)$ is proportional to $\ap$ and 
$e(\mu_i){\cdot}\xi(-\mu_i)$ is proportional to $\mu_i$. 
Since, by Corollary~\ref{I}, $\ap$ lies in the $\Bbb Q$-span of 
$\{\mu_i\mid i\in \{1\}\cup I \}$, 
we may choose the $e(\mu_i)$'s such that these summands will cancel out.

This completes the proof of the proposition.
\end{proof*}%
\begin{s}{Proposition}   \label{net}
For any $a\in\bbk$, the algebra $\p_{\xi_a}$ cannot be a generic stationary 
subalgebra for the coadjoint representation $(P:\p^*)$.
\end{s}\begin{proof*}
Let us show that there us an $h \in \te$ such that $[h,e_a]=e_a$.
(One and the same element for all $e_a$'s).
Choose any $h\in \te$ subject to the requirement that $\mu_i(h)=1$ for $i\in \{1\}\cup I$.
It then follows from Lemma~\ref{koeff}(3) that \ $-\ap(h)=1$ as well.
Hence $(\theta-\ap-\mu_i)(h)=1$, too. But this exactly means that
$[h,e_a]=e_a$. 

Thus, $[\p,\p_{\xi_a}]\cap \p_{\xi_a}\ne 0$ for each $a$.
By \cite[Cor.\,1.8(i)]{tayu1}, this means that $\p_{\xi_a}$ cannot be a generic 
stationary subalgebra.
\end{proof*}%
\begin{s}{Lemma}   \label{dense}
The set $P{\cdot}L$ is dense in $\p^*$.
\end{s}\begin{proof}
It is a standard exercise in Invariant Theory (cf. 
\cite[Lemma\,1]{alela} and \cite[Theorem\,7.3]{VP}).
We have the natural morphism $\phi: P\times L\to \p^*$, $(p, \xi_a)\mapsto p{\cdot}\xi_a$.
It suffices to prove that the differential of $\phi$ is onto at some point.
As such a point, we take $z=(1_P, \xi_0)$, where $1_P$ is the unit of the group $P$.
Then  $d\phi_z(\p, \bbk \xi(\ap))=\p{\cdot}\xi_0 + \bbk \xi(\ap)$.
Here $\p{\cdot}\xi_0$ is a subspace of codimension one in $\p^*$.
Since $\p{\cdot}\xi_0$ is the annihilator of $\p_{\xi_0}=\bbk e_0$, it follows from
Proposition~\ref{stab} that the line $\bbk \xi(\ap)$ is not contained in
$\p{\cdot}\xi_0$. 
\end{proof}%
Now, combining Proposition~\ref{net} and Lemma~\ref{dense}, we complete the proof
of Theorem~\ref{main6} in case, where $\ind\be=0$.

{\bf II.}  \ In general, the argument does not essentially change. 
Now, we have the vector space $\h\subset\te$ of dimension $d:=\rk\g-l$,
and we set 
\[
L=\{\xi_0+x+ a\xi(\ap)\mid x\in\h,\ a\in\bbk\}\subset \p^* \ . 
\]
It is an affine space of dimension $d+1$. By Proposition~\ref{borel},
the projection $\p^* \to \be^*$ takes all points of $L$ to generic
points of $\be^*$. Hence, for any $\xi\in L$ we have $\dim\p{\cdot}\xi\ge 
\dim\be-\dim\h$, i.e., $\dim\p_\xi\le \dim\h+1$. 
On the other hand, we have 

\begin{s}{Proposition}
For any $\xi=\xi_0+x+a\xi(\ap)\in L$, 
we have $\p_\xi \supset \h \oplus\bbk e_a$.
\end{s}%
The proof of Proposition~\ref{stab} goes through verbatim in this
situation, since all the roots involved
are orthogonal to $\h$. 
\\
Thus, we actually have an equality in the last proposition.
Then we prove in the same fashion that $P{\cdot}L$ is dense in $\p^*$ and
neither of $\p_{\xi}$, $\xi\in L$, can be a generic stabiliser.

Thus, Theorem~\ref{main6} is proved. \hfill $\square$

\end{document}